\documentclass[a4paper,11pt,twoside]{article}
\usepackage{isolatin1}
\usepackage{amsfonts,amssymb,amsmath,amsthm}
\usepackage{epic,eepic}
\usepackage{graphicx}
\usepackage[all]{xy}

\DeclareMathOperator{\nr}{nr}
\DeclareMathOperator{\tr}{tr}
\DeclareMathOperator{\Aut}{Aut}

\DeclareMathOperator{\Pic}{Pic}

\newcommand{\R}{{\mathbb{R}}}
\newcommand{\Q}{{\mathbb{Q}}}
\newcommand{\Z}{{\mathbb{Z}}}

\newcommand{\HxH}{{{\cal H}\! \times\! {\cal H}}}
\newcommand{\legendre}[2]{{\left(\frac{#1}{#2}\right)}}
\newcommand{\I}{{\text{I}}}
\newcommand{\II}{{\text{II}}}
\newcommand{\III}{{\text{III}}}

\newcommand{\klammer}{{\raisebox{-1.45ex}[0pt][0pt]{%
      \makebox[0pt][l]{$\Bigl.\ \,\Bigr\}$}}}}
\newcommand{\Klammer}{{\raisebox{-2.7ex}[0pt][0pt]{%
      \makebox[0pt][l]{$\Biggl.\ \,\Biggr\}$}}}}

\theoremstyle{plain}
    \newtheorem{theorem}{Theorem}
    \newtheorem{proposition}[theorem]{Proposition}
    \newtheorem{lemma}[theorem]{Lemma}
    
\theoremstyle{definition}

\begin{document}

\title{On Compact Shimura Surfaces}
\author{Håkan Granath\\
Max-Planck-Institut für Mathematik,\\
Vivatsgasse 7, D-53111 Bonn, Germany\\
email: granath@mpim-bonn.mpg.de}
\date{}
\maketitle

\begin{abstract}
  \noindent
  We study surfaces constructed from groups of units in quaternion
  orders $\Lambda$ over the integers in real quadratic fields~$k$. A
  short presentation of some general theory of such surfaces is given,
  in particular, we construct certain ``modular curves'' on the
  surfaces and examine their properties. We apply our results to some
  surfaces related to the case $k=\Q(\sqrt{13})$ and $d(\Lambda)=(3)$,
  and find their place in the Kodaira classification of algebraic
  surfaces.
  \\
  \hbox{}
  \\
  \textbf{Key words:} Shimura surface, quaternion order, Kodaira
  classification
  \\
  \textbf{2000 Mathematics Subject Classification:}
  14G35, % Modular and Shimura varieties
  11G18. % Arithmetic aspects of modular and Shimura varieties
\end{abstract}

\noindent
The purpose of this paper is to describe some tools to study
particular examples of quaternionic Shimura surfaces from the point of
view of algebraic geometry. To do that one needs of course standard
numerical invariants, but also more explicit knowledge of the geometry
of the curves on these surfaces. In section~\ref{sec:modular-curves},
we construct a family of ``modular curves'', and state some results
regarding these. Proofs of all claims can be found in~\cite{granath}.
In section~\ref{sec:example} we consider an explicit example over the
base field $\Q(\sqrt{13})$. For more details about the underlying
computations, we refer to~\cite{granath}.

%%%%%%%%%%%%%%%%%%%%%%%%%%%%%%%%%%%%%%%%%%%%%%%%%%%%%%%%%%%%%%%%%%%%%%
\section{Introduction}
\label{sec:introduction}

Let $k=\Q(\sqrt{d})$ be a real quadratic field with $R={\cal{O}}_{k}$
and $D=d(k)$. Let $A$ be a totally indefinite quaternion algebra
over~$k$ which allows an involution of type~$2$, i.e.\ a map
$A\rightarrow A$, denoted $x\mapsto\overline{x}$, such that
$\overline{\overline{x}}=x$,
$\overline{x+y}=\overline{x}+\overline{y}$ and
$\overline{x y}=\overline{x}\,\overline{y}$ for all $x,y\in A$, and
such that the restriction of the involution to $k$ is the non-trivial
automorphism of~$k$. If $\Lambda$ is a maximal $R$-order in $A$, then
we let $\Lambda_{0}=\{\lambda\in\Lambda\mid
\overline{\lambda}=\lambda\}$. It can be shown that, replacing the
involution of type 2 with another one if necessary, one can assume
that $D$ and $d(\Lambda_{0})$ are coprime. In that case, we say that
the involution is special with respect to~$\Lambda$.

The canonical involution of a quaternion algebra $A$ is denoted
$x\mapsto x^{*}$. The reduced norm is $\nr(x)=x x^{*}$ and the reduced
trace is $\tr(x)=x+x^{*}$. The reduced discriminant of an order
$\Lambda$ is denoted $d(\Lambda)$, and we let $d_{0}(\Lambda)$ be the
positive rational integer that generates $d(\Lambda)$ if such exists.

Choose a maximal order $\Lambda$ and let
$\Lambda^{1}=\{\lambda\in\Lambda\mid\nr(\lambda)=1\}$. Choose a real
representation $A\rightarrow M_{2}(\R)$, and consider the induced
action of $\Lambda^{1}$ on $\mathcal{H}$ given by Möbius
transformations. We get an action on $\HxH$ by
$\lambda(z_{1},z_{2})=(\lambda z_{1},\overline{\lambda}z_{2})$. We let
$\Gamma$ denote the image of $\Lambda^{1}$ in $\Aut(\HxH)$ and
$X=\HxH/\Gamma$.

$X$ is a so called Hilbert modular surface if $A\cong M_{2}(k)$. These
are not compact, but can be compactified by adding a finite number of
points, the so called cusps. These points give complicated
singularities which were first resolved by Hirzebruch,
see~\cite{hirzebruch5}. Studying numerical invariants of the surfaces
and configurations of certain curves on them, curves coming from the
resolution of the singularities and so called modular curves, one
managed to identify many of these surfaces in terms of the Kodaira
classification of algebraic surface. For example, if
$\Lambda=M_{2}(R)$, then $X$ is rational if and only of $D=5$, 8, 12,
13, 17, 21, 24, 28, 33, or 60 (see~\cite{hirzebruch4}). See
also~\cite{geer1}.

From now on we will consider non-split algebras, i.e.\ we assume that
${A\not\cong M_{2}(k)}$. In this case the quotient
$X=\HxH/\Lambda^{1}$ is compact. However, in general $X$ has quotient
singularities coming from elements of finite order in $\Gamma$, so
called elliptic elements. We let $Y$ denote the canonical minimal
resolution of~$X$. We now summarise some known results about the
numerical invariants of $X$ and $Y$ in this case. Let $\omega$ be the
Gauss-Bonnet form on $\HxH$:
\begin{displaymath}
  \omega=\frac{1}{(2\pi)^{2}}
      \frac{dx_{1}\wedge dy_{1}}{y_{1}^{2}}\wedge
      \frac{dx_{2}\wedge dy_{2}}{y_{2}^{2}},
\end{displaymath}
where $z_{k}=x_{k}+i y_{k}$, $k=1,2$ are the standard coordinates on
the two factors of $\HxH$. The Euler characteristic $e(X)$ is given by
\begin{equation}
  \label{eq:invariants-Euler}
  e(X)=\int_{{\cal{F}}}\omega+\sum_{r>1}e_{r}\frac{r-1}{r},
\end{equation}
where ${\cal{F}}$ is a fundamental domain of the group action and
$e_{r}$ is the number of classes of elliptic points of order $r$
(see~\cite{hirzebruch1}, p.~197). The formula for the volume of a
fundamental domain in $\HxH$ under the action of $\Gamma$ can be found
for example in~\cite{vigneras2}, p.~193. In our case, when $\Lambda$
is a maximal order and the quaternion algebra $A$ admits an involution
of type~$2$, the formula may be written as
\begin{equation}
  \label{eq:invariants-volume}
  \int_{{\cal F}}\omega=2\zeta_{k}(-1)\prod_{p\mid d_{0}(\Lambda)}(p-1)^{2}.
\end{equation}
Here $\zeta_{k}$ is the zeta-function of the field~$k$. The value
$\zeta_{k}(-1)$ can be easily computed using a formula of Siegel (see
e.g.~\cite{hirzebruch1}, p.~192). By Satz~8 in~\cite{freitag}, we have
that the irregularity $q$ of the surface $Y$ vanishes: $q(Y)=0$. We
also have, since $A$ is a skew field, that $\chi(X)=\frac{1}{4}e(X)$
(see~\cite{vigneras2}, p.~206--207).

Apart from the general facts above, not very much is known about the
geometry of specific examples. In~\cite{shavel} all cases were
determined when $X$ is non-singular and has geometric genus $p_{g}=0$.
There are $3$ such cases, one with $k=\Q(\sqrt{2})$ and two with
$k=\Q(\sqrt{3})$. One can also find there results about the
corresponding question for certain extensions of the group~$\Gamma$.
This gives many more cases with $p_{g}=0$. Note however, that in none
of these cases the algebra allows an involution of type 2.

Notation. For $p$ a prime, we let $v_{p}$ denote the valuation on
$\Q_{p}$ and $(a,b)_{p}$ the Hilbert symbol, $a,b\in\Q_{p}^{*}$. For
integers $a,b,c$, we let $\varphi=[a,b,c]$ denote the quadratic form
$a x^{2}+b x y +c y^{2}$. The discriminant of $\varphi$ is
$d(\varphi)=b^{2}-4a c$, and the content $m(\varphi)$ of $\varphi$ is
the ideal generated by $a$, $b$ and~$c$.

%%%%%%%%%%%%%%%%%%%%%%%%%%%%%%%%%%%%%%%%%%%%%%%%%%%%%%%%%%%%%%%%%%%%%%
\section{The modular curves}
\label{sec:modular-curves}

Consider the following 4-dimensional $\Q$-vector space:
\begin{displaymath}
  W=\{\beta\in A\mid \overline{\beta}^{*}=\beta\}.
\end{displaymath}
Clearly, if $\beta\in W$, then $\nr(\beta)\in\Q$. If $\beta\in W$ and
$\nr(\beta)>0$, then we define the following curve on $\HxH$:
\begin{displaymath}
  C_{\beta}=\{(z,\beta z)\mid z\in{\cal H} \}.
\end{displaymath}
This curve maps to a subvariety $F_{\beta}\subset X$. 

We let $\Gamma_{\beta}=\{\lambda\in\Lambda^{1}\mid\lambda\cdot
C_{\beta}=C_{\beta}\}=\{\lambda\in\lambda^{1}\mid
\beta\lambda=\pm\overline{\lambda}\beta\}$, the stabiliser group of
$C_{\beta}$. We define $A_{\beta}=\{a\in A\mid\beta
a=\overline{a}\beta \}$, which is a quaternion algebra over $\Q$, and
let $\Lambda_{\beta}=A_{\beta}\cap\Lambda$. We have that
$\Lambda_{\beta}^{1}$ is a subgroup of $\Gamma_{\beta}$ of index 1
or~2.

If we consider $A$ as a right $A$-module, then we say that a map
$\Phi:A\times A\rightarrow A$ is an $A$-hermitian form if
$\Phi(x+y,z)= \Phi(x,z)+\Phi(y,z)$, $\Phi(x a,y)= \Phi(x,y)a$ and
$\Phi(x,y)= \overline{\Phi(y,x)}^{*}$ for all $a,x,y,z\in A$. If
$\Lambda$ is a maximal order in $A$ and the type $2$ involution on $A$
is special with respect to $\Lambda$, then we say that the
$A$-hermitian form $\Phi$ is integral with respect to the order
$\Lambda$ if
\begin{enumerate}
\item[i)] $\Phi(x,y)\in\Lambda$ for all $x,y\in\Lambda$,
\item[ii)] $\Phi(x,x)\in R+\sqrt{D}\Lambda$ for all $x\in\Lambda$.
\end{enumerate}
For any $\beta$ in $W$, we define a hermitian form $\Phi_{\beta}$ by
$\Phi_{\beta}(x,y)=\overline{y}^{*}\beta x$. We define a $\Z$-lattice
$L$ of rank $4$ by
\begin{equation}
  \label{eq:L}
  L=\{\beta\in W\mid\Phi_{\beta}~\text{is integral}\}.
\end{equation}
Consider a quadratic form $q$ on $L$ by $\displaystyle
q(\beta)=\frac{d_{0}(\Lambda)}{d_{0}(\Lambda)}\nr(\beta)$. We define a
lattice $L^{\#}$ dual to $L$ by
\begin{equation}
  \label{eq:L-dual}
  L^{\#}=\{l\in W\mid\tr(l^{*}L)\subseteq\Z \},
\end{equation}
and a quadratic form $q^{\#}$ on $L^{\#}$ given by $q^{\#}(l)=D
d_{0}(\Lambda_{0})\nr(l)$. One can show that $(L,q)$ and
$(L^{\#},q^{\#})$ are primitive integral quadratic forms, and that
their genera are independent of the choice of involution.

Consider now the quaternary quadratic lattice $(L^{\#},q^{\#})$.
Recall the definition of the even Clifford algebra
$C_{0}(L^{\#},q^{\#})$ as the quotient of the tensor algebra
${\cal{T}}_{0}(L^{\#})=\bigoplus_{k=0}^{\infty}(L^{\#})^{\otimes 2k}$
by the ideal generated by all elements $x\otimes x-q^{\#}(x)$, where
$x\in L^{\#}$. Define a function
$\phi:L^{\#}\otimes_{\Z}L^{\#}\rightarrow A$ by
\begin{displaymath}
  \phi(l_{1}\otimes l_{2})=D d_{0}(\Lambda_{0})l_{1}^{*}l_{2}.
\end{displaymath}
We can extend $\phi$ in a natural way to the even tensor algebra, so
we get a map $\phi:{\cal{T}}_{0}(L^{\#})\rightarrow A$. This map
clearly vanishes on elements of the form $x\otimes x-q^{\#}(x)$, where
$x\in L^{\#}$, giving an embedding of the ring $C_{0}(L^{\#},q^{\#})$
into the algebra $A$. In fact, the image is a subring (but not a
suborder) of~$\Lambda$:
\begin{equation}
  \label{eq:phi}
  \phi:C_{0}(L^{\#},q^{\#})\rightarrow\Lambda.
\end{equation}
The map~\eqref{eq:phi} is one of the main ingredients in the proofs of
theorem~\ref{thm:Lambda-beta-disc} and theorem~\ref{thm:s-phi} below.

Combining the following two results, the genus of the orders
$\Lambda_{\beta}$ is determined.
\begin{proposition}
  \label{prop:S-primitive-orders}
  (Cf.\ prop. 1.12 in~\cite{brzezinski-automorphisms}) Let $F$ be a
  local field, $K$ a separable maximal commutative subalgebra of $A$
  and $S$ a maximal order of~$K$. Let $\Lambda_{1},\Lambda_{2}\subset
  A$ be two $S$-primitive orders (i.e.\ they allow an embedding of
  $S$). If $d(\Lambda_{1})=d(\Lambda_{2})$, then
  $\Lambda_{1}\cong\Lambda_{2}$.
\end{proposition}

\begin{theorem}
  \label{thm:Lambda-beta-disc}
  (See~\cite{granath}, thm. 7.16 and prop. 7.23) The discriminant of
  $\Lambda_{\beta}$ equals the least common multiple of $q(\beta)$ and
  $d_{0}(\Lambda)$. Furthermore, if $p$ is a prime, then the following
  holds:
  \begin{enumerate}
  \item[i\,$)$] If $p$ is split in $k$, then $(A_{\beta})_{p}$ is
    split if and only if $p\nmid d(\Lambda)$. If $p$ is unramified in
    $k$, then $(A_{\beta})_{p}$ is split if and only if
    $v_{p}(q(\beta))$ is even. If $p$ is ramified in $k$, then
    $(A_{\beta})_{p}$ is split if and only if
    $\nr(\beta)\in\nr_{k_{p}/\Q_{p}}(k_{p}^{*})$.
  \item[ii\,$)$] The order $(\Lambda_{\beta})_{p}$ is
    $R_{p}$-primitive unless $p$ is split in $k$ and $p\mid
    d(\Lambda)$.
  \end{enumerate}
\end{theorem}

An element $\beta\in L$ is called primitive if $\{x\in\Q\mid x\beta\in
L\}=\Z$. For $N$ a positive integer, let $L(N)$ denote the set of
primitive elements of $L$ with $q(\beta)=N$.
\begin{proposition}
  \label{prop:group-index}
  (See~\cite{granath}, prop. 8.8) If $\beta\in L(N)$, then
  $[\Gamma_{\beta}:\Lambda_{\beta}^{1}]=2$ if and only if the
  following conditions hold for every prime $p$ such that $p\mid D$
  $:$
  \begin{enumerate}
  \item [i\,$)$] if $(-1,D)_{p}=1$, then $p\mid N$,
  \item [ii\,$)$] if $(-1,D)_{p}=-1$, then $v_{p}(d)\leq
    v_{p}(N)<2v_{p}(D)$.
  \end{enumerate}
\end{proposition}
By proposition~\ref{prop:group-index}, there exists a function
$\kappa$ defined on the set of all positive integers primitively
represented by $(L,q)$, such that
$[\Gamma_{\beta}:\Lambda_{\beta}^{1}]=\kappa(N)$ for all $\beta\in
L(N)$. For every positive integer $N$, we define
\begin{displaymath}
  F_{N}=\bigcup_{\beta\in L(N)}F_{\beta}.
\end{displaymath}
Let $f_{N}$ denote the number of irreducible components of~$F_{N}$.
For every prime $p$ dividing $D$, we define a character
$\chi_{D,p}:\Z_{p}\rightarrow \{\pm1,0\}$, by
$\chi_{D,p}(N)=(D,N)_{p}$ if $p\nmid N$ and $\chi_{D,p}(N)=0$ if
$p\mid N$. Furthermore, if $N\in\Z_{p}$ and $N\neq 0$, then we define
$a_{D,p}(N)=2$ if $v_{p}(N)\geq2v_{p}(D)$ and $a_{D,p}(N)=1$
otherwise.
\begin{theorem}
  \label{thm:f-N}
  (See~\cite{granath}, thm. 8.9) If there exists a prime $p$ such that
  $p\mid d(\Lambda)$ and $p^{2}\mid N$, then $f_{N}=0$. Otherwise
  \begin{displaymath}
    f_{N}=\frac{\kappa(N)}{2}\prod_{p\mid D}
    (\chi_{D,p}(N B)+1)a_{D,p}(N),
  \end{displaymath}
  where $B=d_{0}(\Lambda_{0})/d_{0}(\Lambda)$.
\end{theorem}

Now we have the tools to determine the number of components of $F_{N}$
and their genus. We turn to the problem of determining how these
curves intersect. If $z=(z_{1},z_{2})$ is a point in $\HxH$, then we
let
\begin{displaymath}
  L_{z}=\{\beta\in L\mid q(\beta)>0\text{~and~}
  \beta z_{1}=z_{2}\}\cup\{0\}.
\end{displaymath}
We have that $q$ restricted to $L_{z}$ is positive definite, so the
rank of $L_{z}$ is at most 2. The point $z$ is called special if the
lattice $L_{z}$ has rank~2. It is possible to choose an orientation on
the 2-dimensional $\R$-vector bundle ${\cal{L}}=\{(z,\beta)\in
(\HxH)\times M_{2}(\R)\mid {\det(\beta)>0}\text{~and~}\beta
z_{1}=z_{2}\}\cup \HxH\times\{0\}$ on $\HxH$. We fix such a choice,
and get induced orientations on the lattices~$L_{z}$. Let
$q_{z}=q|_{L_{z}}$ be the oriented binary form associated with the
special point~$z$.

We introduce (as in~\cite{hirzebruch3}) a function $s$ counting the
number of special points. Let $\varphi$ be an oriented binary form.
Then we put
\begin{equation}
  \label{eq:def-of-s-phi}
  s(\varphi)=\sum_{\zeta\in X(\varphi)}
  \frac{w_{\varphi}}{v_{\zeta}},
\end{equation}
where $w_{\varphi}$ is the number of elements in the group of oriented
automorphisms of $\varphi$, and $v_{\zeta}$ is the order of the
isotropy group $\Lambda_{z}^{1}$ in $\Lambda^{1}$ of a point
$z\in\HxH$ representing $\zeta$ (so $w_{\varphi}$ and $v_{\zeta}$ both
take values in the set $\{2,4,6\}$).

From now on, we will assume that $(D,6)=1$. In this case, every
elliptic point is of order 2, 3 or 5 (the latter case is only possible
if $k=\Q(\sqrt{5})$), and one can show that every elliptic point of
order 2 or 3 is a special point.
\begin{proposition}
  \label{prop:binary-form-at-elliptic}
  (See~\cite{granath}, cor. 9.9) Let $q_{z}$ be the quadratic form of
  an elliptic point $z$ of order $n=2$ or $n=3$. There are two cases:
  \begin{enumerate}
  \item[i\,$)$] $q_{z}\cong m D\varphi$, where $\varphi=[1,0,1]$ if
    $n=2$ or $\varphi=[1,-1,1]$ if $n=3$, and where $m$ is a positive
    integer with $m\mid d(\Lambda)$. If $n=2$, then we have
    $\Lambda^{1}_{z}\cap\Lambda_{\beta}^{1}= \{\pm1\}$ and
    $\Lambda^{1}_{z}\subseteq\Gamma_{\beta}$ for all $\beta\in
    L_{z}\setminus\{0\}$, and if $n=3$, then we have
    $\Lambda^{1}_{z}\cap\Gamma_{\beta}=\{\pm1\}$ for all $\beta\in
    L_{z}\setminus\{0\}$.
  \item[ii\,$)$] $q_{z}\cong m\varphi$, where $d_{0}(\varphi)=-4D$ if
    $n=2$ or $d_{0}(\varphi)=-3D$ if $n=3$, and $m$ is a positive
    integer with $m\mid d(\Lambda)$. We have $\Lambda^{1}_{z}\subseteq
    \Lambda_{\beta}^{1}$ for all $\beta\in L_{z}\setminus\{0\}$.
  \end{enumerate}
\end{proposition}

% Tala om denna rättelse för Juliusz
\begin{lemma}
  \label{lem:which_binary-locally}
  (See~\cite{granath}, lem. 9.5) Let $p$ be a prime and $\varphi_{p}$
  a non-degenerate binary quadratic form over~$\Z_{p}$. Then
  $(L_{p},q_{p})$ represents $\varphi_{p}$ primitively if and only if
  the following conditions are satisfied:
  \begin{enumerate}
  \item[i\,$)$] if $p$ is unramified in $k$, then $p\nmid
    m(\varphi_{p})$,
  \item[ii\,$)$] if $p\mid d(\Lambda)$, then $\varphi_{p}$ is
    anisotropic, and it is either a modular form with $p^{2}\nmid
    m(\varphi_{p})$ or a non-modular form such that
    $d(\varphi_{p})=(p)$ if $p\neq 2$ and such that
    $d(\varphi_{p})\mid8$ and $\varphi_{p}$ do not primitively
    represent $4$ if $p=2$,
  \item[iii\,$)$] if $p\mid D$, then $p\mid d(\varphi_{p})$ and
    $p^{2}\nmid m(\varphi_{p})$. Furhermore, if $m(\varphi_{p})=(1)$,
    then $\varphi_{p}$ represents
    $B=d_{0}(\Lambda_{\Z})/d_{0}(\Lambda)$, and if
    $m(\varphi_{p})=(p)$, then $\varphi_{p}$ represents $-B D$.
  \end{enumerate}
\end{lemma}

If $\varphi$ is a positive definite oriented binary form such that
$\varphi\cong q_{z}$ for some point $z\in\HxH$, then we let
$\Delta=d(\varphi)$, and $m>0$ be a generator of $m(\varphi)$. We
write $m=m_{1}m_{2}m_{3}$, where $m_{1}$ contains those prime factors
which divide $D$, $m_{2}$ those prime factors which divide
$d(\Lambda)$ and $m_{3}$ contains all other prime factors. By
lemma~\ref{lem:which_binary-locally}, we have that $m_{1}\mid D$, that
$m_{2}\mid d(\Lambda)$ and that every prime dividing $m_{3}$ is split
in~$k$. If $\varphi$ is a binary form over $\Z$ and $p$ a prime, then
we let $\alpha_{p}(\varphi)=1$ if $\varphi_{p}$ is modular, and
$\alpha_{p}(\varphi)=2$ otherwise. We also introduce the modified
class number $h'$ (compare e.g.~\cite{hirzebruch3}), which counts
primitive binary forms $\varphi$ with a given discriminant $N<0$ with
multiplicity $2/w_{\varphi}$. Hence $h'(-3)=1/3$, $h'(-4)=1/2$ and
$h'(N)=h(N)$ otherwise.
\begin{theorem}
  \label{thm:s-phi}
  (See~\cite{granath}, thm. 9.16) A positive definite oriented binary
  form $\varphi$ is represented by $q$ if and only if $\varphi_{p}$
  satisfies the conditions of lemma~\ref{lem:which_binary-locally} for
  every prime~$p$. For such a form, we have
  \begin{equation}
    \label{eq:s-phi}
    s(\varphi)=2^{a-1}h'\left(\frac{\Delta}{m_{2}^{2}D}\right)
    \prod_{p\mid D}\alpha_{p}(\varphi)
    \prod_{p\mid d(\Lambda)}\frac{2}{\alpha_{p}(\varphi)},
  \end{equation}
  where $a$ is the number of different primes dividing~$m_{3}$.
\end{theorem}

Theorem~\ref{thm:s-phi} is our tool to determine intersection points
of the modular curves. It remains to determine the self-intersections
of the curves $F_{\beta}$. For that one can use the theory of Chern
divisors (see~\cite{hirzebruch1}, p.~252).

%%%%%%%%%%%%%%%%%%%%%%%%%%%%%%%%%%%%%%%%%%%%%%%%%%%%%%%%%%%%%%%%%%%%%%
\section{Some surfaces over $\Q(\sqrt{13})$}
\label{sec:example}

We will now consider the example where $k=\Q(\sqrt{13})$ and
$d_{0}(\Lambda)=3$. This is one of the cases with relatively small
hyperbolic volume of the fundamental domain. Since $h(k)=1$ in this
case, there is only one maximal order in $A$ up to isomorphism (by
Eichler's norm theorem and lemma~III.5.6 in~\cite{vigneras1}).

We consider first discrete groups of $\Aut(\HxH)$ extending~$\Gamma$.
Note that every element in $A^{+}=\{a\in A\mid\nr(a)\text{~is totally
  positive}\}$ acts on $\HxH$, giving a map
$\varrho:A^{+}\rightarrow\Aut(\HxH)$. The totally positive elements of
the normaliser of $\Lambda$,
\begin{displaymath}
  N^{+}(\Lambda)=\{n\in A^{+}\mid n\Lambda n^{-1}=\Lambda\},
\end{displaymath}
naturally act on~$X$. Let $\Gamma_{\II}=\varrho(N^{+}(\Lambda))$. In
our case, we have that $\Gamma_{\II}$ is given as follows: Let
$v=4+\sqrt{13}$, so $3=v\overline{v}$. There exists
$\lambda_{x}\in\Lambda$ such that $\nr(\lambda_{x})=x$ for $x=v$,
$\overline{v}$ and $3$, and we let $\gamma_{x}=
\varrho(\lambda_{x}):\HxH\rightarrow\HxH$. This gives three well
defined involutions on $X$, which we denote $\iota_{3}$, $\iota_{v}$
and~$\iota_{\overline{v}}$. We define $\Gamma_{\I}=
\Gamma\cup\gamma_{3}\Gamma$, and we have in fact that $\Gamma_{\II}=
\Gamma_{\I}\cup\gamma_{v}\Gamma_{\I}$.

We will now define an extension of $\Gamma_{\II}$. Let us fix an
element $\beta_{13}\in L$ with $q(\beta_{13})=13$, and define a map
$T_{0}:\HxH\rightarrow\HxH$, by
\begin{displaymath}
  T_{0}(z_{1},z_{2})=(\beta_{13}^{*} z_{2},\beta_{13} z_{1}).
\end{displaymath}
We have that $T_{0}^{2}$ is the identity, and one can check that
$T_{0}\Gamma=\Gamma T_{0}$, so $T_{0}$ induces an involution
on $X$, which we denote~$T$. We have
$T\iota_{x}=\iota_{\overline{x}}T$, for $x\in\{1,v,\overline{v},3\}$.
Let $\Gamma_{\III}=\Gamma_{\II}\cup T_{0}\Gamma_{\II}$. This group
is independent of the choice of the element~$\beta_{13}$. Hence, we
have a tower of discrete groups acting on $\HxH$
\begin{equation}
  \label{eq:Gamma-tower}
  \Gamma\subset\Gamma_{\I}\subset\Gamma_{\II}\subset\Gamma_{\III}.
\end{equation}
We will consider the $4$ quotient surfaces $X$,
$X_{\I}=\HxH/\Gamma_{\I}$, $X_{\II}=\HxH/\Gamma_{\II}$ and
$X_{\III}=\HxH/\Gamma_{\III}$. We let $Y$, $Y_{\I}$, $Y_{\II}$ and
$Y_{\III}$ respectively denote the canonical minimal resolution of the
corresponding quotient surface. We let
$G_{\I}=\{\iota_{1},\iota_{3}\}$,
$G_{\II}=\{\iota_{1},\iota_{v},\iota_{\overline{v}},\iota_{3}\}$ and
$G_{\III}=G_{\II}\cup TG_{\II}$ denote the corresponding groups acting
on~$X$, so $G_{\III}\cong D_{4}$. We have that $G_{\I}$ is the center
of $G_{\III}$, so the filtration~(\ref{eq:Gamma-tower}) is in some
sense canonical. Our main result is:
\begin{theorem}
  \label{thm:classification-of-the-surfaces}
  (See~\cite{granath}, thm. 12.1) $Y$ is a minimal surface of general
  type, $Y_{\I}$ is a $K3$-surface blown up $4$ times, $Y_{\II}$ is an
  Enriques surface blown up $2$ times and $Y_{\III}$ is a rational
  surface with Euler number $e=12$.
\end{theorem}

Calculation of the quadratic form $q$ in this case gives:
\begin{equation}
  \label{eq:example-q}
  q = 2 {t_{0}}^{2}+t_{0} t_{1}+5 {t_{1}}^{2}
     -13({t_{2}}^{2} - t_{2} t_{3} +{t_{3}}^{2})
\end{equation}
By theorem~\ref{thm:f-N}, the number of components of $F_{N}$ is
\begin{equation}
  \label{eq:example-f-N}
  f_{N}=
  \begin{cases}
    0 & \text{if $\legendre{N}{13}=1$ or $9\mid N$,} \\
    2 & \text{if $13^{2}\mid N$ and $9\nmid N$,} \\
    1 & \text{otherwise.}
  \end{cases}
\end{equation}
If $\beta$ is a primitive element of $L$, then we have, by
theorem~\ref{thm:Lambda-beta-disc}, that $d_{0}(\Lambda_{\beta})$ is
the least common multiple of $q(\beta)$ and~$3$. By
proposition~\ref{prop:group-index}, we have
\begin{displaymath}
  \kappa(N)=
  \begin{cases}
    2 & \text{if $13\mid N$,} \\
    1 & \text{otherwise.}
  \end{cases}
\end{displaymath}
The action of $G_{\III}$ on the curves are given by
\begin{lemma}
  (See~\cite{granath}, lem. 12.4) For every positive integer $N$
  represented by $q$, we have that $\iota_{3}(F_{N})=F_{N}$ and
  $T(F_{N})=F_{N}$.  If $3\nmid N$, then $\iota_{x}(F_{N})=F_{3N}$\ 
  and $\iota_{x}(F_{3N})=F_{N},$ for $x=v$ and $x=\overline{v}$.
\end{lemma}

Now, $X$ has 8 elliptic points of order 2, and 4 of order 3,
corresponding to elliptic elements $\epsilon_{1},\dots,\epsilon_{8}$
and $\rho_{1},\dots,\rho_{4}$ respectively. The elements of order 2
give rise to $(-2)$-curves on $Y$, $\rho_{1},\rho_{2}$ to two
intersecting $(-2)$-curves each and $\rho_{3},\rho_{4}$ to
$(-3)$-curves. In addition to these curves, we will consider the 6
irreducible curves $F_{2}$, $F_{6}$, $F_{5}$, $F_{15}$, $F_{13}$
and~$F_{39}$. The binary forms associated to the elliptic points, and
the branches of the above curves that pass through them, are given by
table~\ref{tab:elliptic-points}. All other intersection points of
these curves are given in table~\ref{table:special-points}. We get
that $F_{2}$ and $F_{6}$ are non-singular rational curves, $F_{13}$
and $F_{39}$ are non-singular elliptic curves and $F_{5}$ and $F_{15}$
are elliptic curves with two nodes each.

Determining the fixed points of the elements of $G_{\III}$ requires
explicit calculations. The isotropy groups of all points with
non-trivial isotropy group $\Gamma_{\II,z}$ are given in
table~$\ref{tab:non-trivial-isotropy}$. $F_{13}$ and $F_{39}$ are
pointwise fixed by $T$ and $T\iota_{3}$ respectively. All other fixed
points of elements in $G_{\II}T$ are already fixed points of elements
of $G_{\II}$.
\begin{table}
  \begin{center}
    \renewcommand{\arraystretch}{1.1}
    \begin{tabular}{cccc}
      Elliptic element & Order & Binary form & Branches \\
      \hline
      $\epsilon_{1},\ \epsilon_{2}$ & 2 & $[13,0,13]$ & $F_{13}$, $F_{13}$ \\
      $\epsilon_{3},\ \epsilon_{4}$ & 2 & $[39,0,39]$ & $F_{39}$, $F_{39}$ \\
      $\epsilon_{5},\ \epsilon_{6}$ & 2 & $[2,2,7]$ & $F_{2}$ \\
      $\epsilon_{7},\ \epsilon_{8}$ & 2 & $[6,6,21]$ & $F_{6}$ \\
      $\rho_{1},\ \rho_{2}$ & 3 & $[13,-13,13]$ & $F_{13}$, $F_{39}$ \\
      $\rho_{3}$ & 3 & $[2,1,5]$ & $F_{2}$, $F_{6}$, $F_{5}$, $F_{15}$  \\
      $\rho_{4}$ & 3 & $[2,-1,5]$ & $F_{2}$, $F_{6}$, $F_{5}$, $F_{15}$  \\
    \end{tabular}
    \caption{Elliptic points}
    \label{tab:elliptic-points}
  \end{center}
\end{table}

\begin{table}
  \begin{center}
    \renewcommand{\arraystretch}{1.1}
    \begin{tabular}{ccc}
      Binary form & Branches & $s(\varphi)$ \\
      \hline
      $[2, 0, 39]$ & $F_{2}$, $F_{39}$ & $h(-24)=2$ \\
      $[6, 0, 13]$ & $F_{6}$, $F_{13}$ & $h(-24)=2$ \\
      $[5, -3, 5]$ & $F_{5}$, $F_{5}$, $F_{13}$ & $2h(-7)=2$ \\
      $[15, -9, 15]$ & $F_{15}$, $F_{15}$, $F_{39}$ & $2h(-7)=2$ \\
      $[5, 2, 8]$ & $F_{5}$, $F_{15}$ & $h(-12)=1$ \\
      $[5, -2, 8]$ & $F_{5}$, $F_{15}$ & $h(-12)=1$ \\
      $[5, 0, 39]$ & $F_{5}$, $F_{39}$ & $h(-60)=2$ \\
      $[15, 0, 13]$ & $F_{15}$, $F_{13}$ & $h(-60)=2$ \\
      $[5, -5, 11]$ & $F_{5}$, $F_{39}$ & $h(-15)=2$ \\
      $[7, -1, 7]$ & $F_{15}$, $F_{13}$ & $h(-15)=2$ \\
      $[13, 0, 39]$ & $F_{13}$, $F_{39}$ & $h(-156)/2=2$ \\
    \end{tabular}
    \caption{Intersection points}
    \label{table:special-points}
  \end{center}
\end{table}

\begin{table}[htbp]
  \centering
  \begin{tabular}{c|@{\ \ \ }c@{\ \ \ \ \ }c@{\ \ \ \ \ }c@{\ \ \ \ \ }c}
    \raisebox{0ex}[0ex][1.0ex]{}
    & $\Gamma_{z}$ & $\Gamma_{\I,z}$ & $\Gamma_{\II,z}$ & $\Gamma_{\III,z}$ \\
    \hline
    \raisebox{0ex}[2.7ex]{}%
      $[2,-2,7]$ & $\Z/2$\klammer &
      \raisebox{-1.6ex}[0pt][0pt]{$\Z/2$\Klammer} & & \\
    $[2,-2,7]$ & $\Z/2$ & & \raisebox{-1.6ex}[0pt][0pt]{$\Z/2$}  &
      \raisebox{-1.6ex}[0pt][0pt]{$\Z/4$} \\
    $[6,-6,21]$ & $\Z/2$\klammer & \raisebox{-1.6ex}[0pt][0pt]{$\Z/2$} & &  \\
    \raisebox{0ex}[0ex][1.5ex]{}%
      $[6,-6,21]$ & $\Z/2$ &  &  &  \\
    \hline 
    \raisebox{0ex}[2.7ex]{}%
      $[13,0,13]$ & $\Z/2$\klammer &
      \raisebox{-1.6ex}[0pt][0pt]{$\Z/2$\Klammer} & & \\
    $[13,0,13]$ & $\Z/2$ & & \raisebox{-1.6ex}[0pt][0pt]{$\Z/2$} &
      \raisebox{-1.6ex}[0pt][0pt]{$\Z/2\times\Z/2$} \\
    $[39,0,39]$ & $\Z/2$\klammer & \raisebox{-1.6ex}[0pt][0pt]{$\Z/2$} & &  \\
    \raisebox{0ex}[0ex][1.5ex]{}%
      $[39,0,39]$ & $\Z/2$ &  &  &  \\
    \hline
    \raisebox{0ex}[2.7ex]{}%
      $[2,1,5]$ & $\Z/3$ & $\Z/6$\klammer &
      \raisebox{-1.6ex}[0pt][0pt]{$\Z/6$} &
      \raisebox{-1.6ex}[0pt][0pt]{$\Z/12$} \\
    \raisebox{0ex}[0ex][1.5ex]{}%
      $[2,-1,5]$ & $\Z/3$ & $\Z/6$ &  &  \\
    \hline
    \raisebox{0ex}[2.7ex]{}%
      $[13,-13,13]$ & $\Z/3$ & $\Z/6$\klammer &
      \raisebox{-1.6ex}[0pt][0pt]{$\Z/6$} &
      \raisebox{-1.6ex}[0pt][0pt]{$D_{6}$} \\
    \raisebox{0ex}[0ex][1.5ex]{}%
    $  [13,-13,13]$ & $\Z/3$ & $\Z/6$ &  &  \\
    \hline
    \raisebox{0ex}[2.7ex]{}%
      $[13,0,39]$ & $(1)$ & $\Z/2$\klammer &
      \raisebox{-1.6ex}[0pt][0pt]{$\Z/2$} &
      \raisebox{-1.6ex}[0pt][0pt]{$\Z/2\times\Z/2$} \\
    \raisebox{0ex}[0ex][1.5ex]{}%
      $[13,0,39]$ & $(1)$ & $\Z/2$ &  &  \\
    \hline
    \raisebox{0ex}[2.7ex]{}%
      $[5,2,8]$ & $(1)$ & $\Z/2$\klammer &
      \raisebox{-1.6ex}[0pt][0pt]{$\Z/2$} &
      \raisebox{-1.6ex}[0pt][0pt]{$\Z/4$} \\
    $[5,-2,8]$ & $(1)$ & $\Z/2$ &  &  \\
  \end{tabular}
  \caption{Overview of points $z$ with non-trivial isotropy group $\Gamma_{\II,z}$}
  \label{tab:non-trivial-isotropy}
\end{table}

We now present the configurations of curves on the different surfaces.
We start with the surface $Y$ and proceed by a sequence of blow ups
and quotient constructions. We will use the convention that whenever
we have two surfaces $Z$ and $\widetilde{Z}$, then $\widetilde{Z}$ is
isomorphic to the surface $Z$ blown up in a finite set of points. For
easier navigation among the surfaces, we first summarise in a diagram
the sequence of surfaces that we will construct:
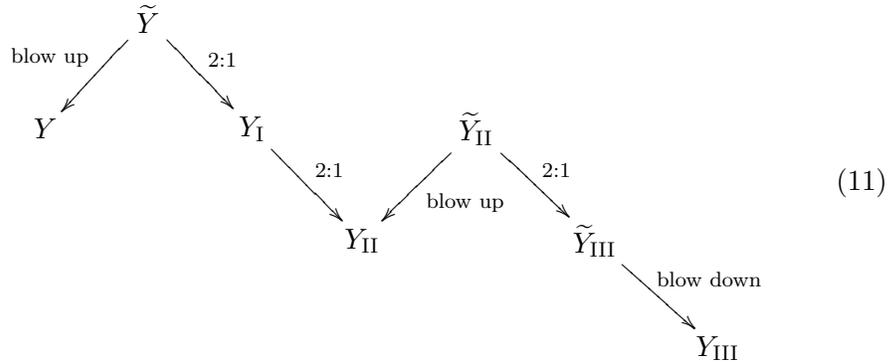
\begin{equation}
  \label{eq:surfaces}
  \begin{split}
    \xymatrix{%
      & \ar[ld]_{\text{blow up}} \widetilde{Y} \ar[rd]^{2:1} \\
      Y & & Y_{\I}\ar[rd]^{2:1} & & \ar[ld]^{\text{blow up}}
      \widetilde{Y}_{\II}\ar[rd]^{2:1} \\ 
      & & & Y_{\II} & &
      \widetilde{Y}_{\III}\ar[rd]^{\text{blow down}} \\
      & & & & & & Y_{\III}
      }
  \end{split}
\end{equation}

We have drawn the configuration of curves on $Y$ in
figure~\ref{fig:Y}. Here all curves drawn with thick lines are
exceptional curves coming from the canonical resolution of the
singularities of~$X$. If the self-intersection number of such a curve
is not explicitly given in the figure, then it is~$-2$. Furthermore,
the rational curves $F_{2}$ and $F_{6}$ have self-intersection $-2$
and the elliptic curves $F_{13}$ and $F_{39}$ have
self-intersection~$-4$. The nodal elliptic curves $F_{5}$ and $F_{15}$
have self-intersection~$2$.

%%%%%%%%%%%%%%%%%%%%%%%%%%%%%%%%%%%%%%%%%%%%%%%%%%%%%%%%%%%%
% Y
\begin{figure}[htbp]
  \centering
  \setlength{\unitlength}{0.9mm}
  \begin{picture}(135,86)(0,-42)
    % F2
    \path(0,25)(70,25)
    \put(35,28){\makebox(0,0)[r]{$F_{2}$}}
    % F6
    \path(0,-25)(70,-25)
    \put(35,-28){\makebox(0,0)[r]{$F_{6}$}}
    % F39
    \put(5,30){\arc{20}{4.2}{6.28}}
    \path(15,30)(15,20)
    \put(25,20){\arc{20}{1.57}{3.14}}
    \path(25,10)(45,10)
    \put(45,20){\arc{20}{0}{1.57}}
    \path(55,20)(55,30)
    \put(65,30){\arc{20}{3.14}{6.28}}
    \path(75,30)(75,10)
    \spline(75,10)(75,5)(78,3)
    \path(78,3)(87,-3)
    \spline(87,-3)(90,-5)(93,-3)
    \path(93,-3)(102,3)
    \spline(102,3)(105,5)(108,3)
    \path(108,3)(117,-3)
    \spline(117,-3)(120,-5)(123,-3)
    \path(123,-3)(135,5)
    \put(16,35){\makebox(0,0)[l]{$F_{39}$}}
    % F13
    \put(5,-30){\arc{20}{0}{2.08}}
    \path(15,-30)(15,-20)
    \put(25,-20){\arc{20}{3.14}{4.71}}
    \path(25,-10)(45,-10)
    \put(45,-20){\arc{20}{4.71}{6.28}}
    \path(55,-20)(55,-30)
    \put(65,-30){\arc{20}{0}{3.14}}
    \path(75,-30)(75,-10)
    \spline(75,-10)(75,-5)(78,-3)
    \path(78,-3)(87,3)
    \spline(87,3)(90,5)(93,3)
    \path(93,3)(102,-3)
    \spline(102,-3)(105,-5)(108,-3)
    \path(108,-3)(117,3)
    \spline(117,3)(120,5)(123,3)
    \path(123,3)(135,-5)
    \put(16,-35){\makebox(0,-0)[l]{$F_{13}$}}
    % F5
    \path(0,-5)(5,-5)
    \spline(5,-5)(10,-5)(13,-1)
    \path(13,-1)(22,11)
    \spline(22,11)(25,15)(28,12)
    \path(28,12)(32,8)
    \spline(32,8)(35,5)(38,8)
    \path(38,8)(42,12)
    \spline(42,12)(45,15)(48,11)
    \path(48,11)(57,-1)
    \spline(57,-1)(60,-5)(64,-5)
    \path(64,-5)(67,-5)
    \spline(67,-5)(70,-5)(73,-8)
    \path(73,-8)(75,-10)
    \path(75,-10)(77.8,-7.49)(79.02,-6.75)(80.06,-6.31)(80.98,-6.08)(81.83,-6)(82.6,-6.04)(83.31,-6.19)(83.95,-6.43)(84.51,-6.74)(85.01,-7.12)(85.43,-7.56)(85.77,-8.04)(86.03,-8.56)(86.21,-9.1)(86.3,-9.66)(86.31,-10.23)(86.23,-10.79)(86.07,-11.33)(85.83,-11.86)(85.5,-12.35)(85.1,-12.79)(84.62,-13.19)(84.06,-13.52)(83.44,-13.77)(82.75,-13.93)(81.99,-14)(81.16,-13.95)(80.25,-13.75)(79.24,-13.36)(78.07,-12.69)(75,-10)
    \path(75,-10)(73,-12)
    \spline(73,-12)(70,-15)(70,-17.5)
    \spline(70,-17.5)(70,-20)(73,-23)
    \path(73,-23)(75,-25)
    \path(75,-25)(77.8,-22.49)(79.02,-21.75)(80.06,-21.31)(80.98,-21.08)(81.83,-21)(82.6,-21.04)(83.31,-21.19)(83.95,-21.43)(84.51,-21.74)(85.01,-22.12)(85.43,-22.56)(85.77,-23.04)(86.03,-23.56)(86.21,-24.1)(86.3,-24.66)(86.31,-25.23)(86.23,-25.79)(86.07,-26.33)(85.83,-26.86)(85.5,-27.35)(85.1,-27.79)(84.62,-28.19)(84.06,-28.52)(83.44,-28.77)(82.75,-28.93)(81.99,-29)(81.16,-28.95)(80.25,-28.75)(79.24,-28.36)(78.07,-27.69)(75,-25)
    \path(75,-25)(71,-29)
    \put(88,-25){\makebox(0,0)[l]{$F_{5}$}}
    % F15
    \path(0,5)(5,5)
    \spline(5,5)(10,5)(13,1)
    \path(13,1)(22,-11)
    \spline(22,-11)(25,-15)(28,-12)
    \path(28,-12)(32,-8)
    \spline(32,-8)(35,-5)(38,-8)
    \path(38,-8)(42,-12)
    \spline(42,-12)(45,-15)(48,-11)
    \path(48,-11)(57,1)
    \spline(57,1)(60,5)(64,5)
    \path(64,5)(67,5)
    \spline(67,5)(70,5)(73,8)
    \path(73,8)(75,10)
    \path(75,10)(77.8,7.49)(79.02,6.75)(80.06,6.31)(80.98,6.08)(81.83,6)(82.6,6.04)(83.31,6.19)(83.95,6.43)(84.51,6.74)(85.01,7.12)(85.43,7.56)(85.77,8.04)(86.03,8.56)(86.21,9.1)(86.3,9.66)(86.31,10.23)(86.23,10.79)(86.07,11.33)(85.83,11.86)(85.5,12.35)(85.1,12.79)(84.62,13.19)(84.06,13.52)(83.44,13.77)(82.75,13.93)(81.99,14)(81.16,13.95)(80.25,13.75)(79.24,13.36)(78.07,12.69)(75,10)
    \path(75,10)(73,12)
    \spline(73,12)(70,15)(70,17.5)
    \spline(70,17.5)(70,20)(73,23)
    \path(73,23)(75,25)
    \path(75,25)(77.8,22.49)(79.02,21.75)(80.06,21.31)(80.98,21.08)(81.83,21)(82.6,21.04)(83.31,21.19)(83.95,21.43)(84.51,21.74)(85.01,22.12)(85.43,22.56)(85.77,23.04)(86.03,23.56)(86.21,24.1)(86.3,24.66)(86.31,25.23)(86.23,25.79)(86.07,26.33)(85.83,26.86)(85.5,27.35)(85.1,27.79)(84.62,28.19)(84.06,28.52)(83.44,28.77)(82.75,28.93)(81.99,29)(81.16,28.95)(80.25,28.75)(79.24,28.36)(78.07,27.69)(75,25)
    \path(75,25)(71,29)
    \put(88,25){\makebox(0,0)[l]{$F_{15}$}}
    % Fixpunkter
    \put(15,0){\makebox(0,0)[l]{$x_{1}$}}
    \put(55,0){\makebox(0,0)[r]{$x_{2}$}}
    \put(85,0){\makebox(0,0)[l]{$x_{3}$}}
    \put(100,0){\makebox(0,0)[l]{$x_{4}$}}
    \put(115,0){\makebox(0,0)[l]{$x_{5}$}}
    \put(130,0){\makebox(0,0)[l]{$x_{8}$}}
    \put(108,10.6){\circle*{.8}}
    \put(109,10.6){\makebox(0,0)[l]{$x_{6}$}}
    \put(108,-10.6){\circle*{.8}}
    \put(109,-10.6){\makebox(0,0)[l]{$x_{7}$}}
    \put(123,10.6){\circle*{.8}}
    \put(124,10.6){\makebox(0,0)[l]{$x_{9}$}}
    \put(123,-10.6){\circle*{.8}}
    \put(124,-10.6){\makebox(0,0)[l]{$x_{10}$}}
    \Thicklines
    % (3;1,1)
    \path(5,-30)(5,30)
    \put(6,10){\makebox(0,0)[l]{$-3$}}
    \put(4,10){\makebox(0,0)[r]{$B_{1}$}}
    \path(65,-30)(65,30)
    \put(64,10){\makebox(0,0)[r]{$-3$}}
    \put(66,10){\makebox(0,0)[l]{$B_{2}$}}
    % (2;1,1)
    \path(10,15)(10,30)
    \put(10,32){\makebox(0,0){$A_{1}$}}
    \path(60,15)(60,30)
    \put(60,32){\makebox(0,0){$A_{2}$}}
    \path(10,-15)(10,-30)
    \put(10,-12){\makebox(0,0){$A_{3}$}}
    \path(60,-15)(60,-30)
    \put(60,-12){\makebox(0,0){$A_{4}$}}
    % (3;1,2)
    \path(107,13)(118,-13)
    \path(107,-13)(118,13)
    \path(122,13)(133,-13)
    \path(122,-13)(133,13)
    % (2;1,1)
    \put(7.5,-45){\arc{15}{3.14}{6.28}}
    \put(65,-45){\arc{15}{3.14}{6.28}}
    \put(7.5,45){\arc{15}{0}{3.14}}
    \put(65,45){\arc{15}{0}{3.14}}
  \end{picture}
  \caption{$Y$}
  \label{fig:Y}
\end{figure}

The involution $\iota_{3}$ has $8$ fixed points on~$X$. On the surface
$Y$, the fixed point locus of the lifted involution consists of the
two $(-3)$-curves and the $10$ isolated points $x_{1},\dots,x_{10}$,
which are indicated in figure~\ref{fig:Y}. Blowing up these 10 points
on $Y$ to the surface~$\widetilde{Y}$, and thereafter taking the
quotient, we get a surface which, in fact, is the minimal
desingularisation $Y_{\I}$ of~$X_{\I}$. We indicate the configuration
of curves in figure~\ref{fig:Y-I}. We have that
$F_{2}^{2}=F_{6}^{2}=-1$, that $F_{13}$ and $F_{39}$ are non-singular
rational curves with $F_{13}^{2}=F_{39}^{2}=-4$, and that $F_{5}$ and
$F_{15}$ are nodal rational curves with $F_{5}^{2}=F_{15}^{2}=0$.

%%%%%%%%%%%%%%%%%%%%%%%%%%%%%%%%%%%%%%%%%%%%%%%%%%%%%%%%%%%%
% Y_{\I}
\begin{figure}[htbp]
  \centering
  \setlength{\unitlength}{0.9mm}
  \begin{picture}(135,86)(0,-42)
    % F2
    \path(0,25)(70,25)
    \put(35,28){\makebox(0,0)[r]{$F_{2}$}}
    % F6
    \path(0,-25)(70,-25)
    \put(35,-28){\makebox(0,0)[r]{$F_{6}$}}
    % F39
    \path(25,-10)(45,-10)
    \put(45,-20){\arc{20}{4.71}{6.28}}
    \path(55,-20)(55,-30)
    \put(65,-30){\arc{20}{0}{3.14}}
    \path(75,-30)(75,-15)
    \put(85,-15){\arc{20}{3.14}{4.71}}
    \path(85,-5)(135,-5)
    \put(54,35){\makebox(0,0)[r]{$F_{39}$}}
    % F13
    \path(25,10)(45,10)
    \put(45,20){\arc{20}{0}{1.57}}
    \path(55,20)(55,30)
    \put(65,30){\arc{20}{3.14}{6.28}}
    \path(75,30)(75,15)
    \put(85,15){\arc{20}{1.57}{3.14}}
    \path(85,5)(135,5)
    \put(54,-35){\makebox(0,0)[r]{$F_{13}$}}
    % F15
    \path(0,5)(27,5)
    \spline(27,5)(30,5)(33,8)
    \path(33,8)(37,12)
    \spline(37,12)(40,15)(43,12)
    \path(43,12)(53.7,1.3) % Glappet!
    \path(56.3,-1.3)(75,-20)
    \path(75,-20)(77.8,-17.49)(79.02,-16.75)(80.06,-16.31)(80.98,-16.08)(81.83,-16)(82.6,-16.04)(83.31,-16.19)(83.95,-16.43)(84.51,-16.74)(85.01,-17.12)(85.43,-17.56)(85.77,-18.04)(86.03,-18.56)(86.21,-19.1)(86.3,-19.66)(86.31,-20.23)(86.23,-20.79)(86.07,-21.33)(85.83,-21.86)(85.5,-22.35)(85.1,-22.79)(84.62,-23.19)(84.06,-23.52)(83.44,-23.77)(82.75,-23.93)(81.99,-24)(81.16,-23.95)(80.25,-23.75)(79.24,-23.36)(78.07,-22.69)(75,-20)
    \path(75,-20)(73,-22)
    \put(33,-5){\makebox(0,0)[l]{$F_{15}$}}
    % % F5
    \path(0,-5)(27,-5)
    \spline(27,-5)(30,-5)(33,-8)
    \path(33,-8)(37,-12)
    \spline(37,-12)(40,-15)(43,-12)
    \path(43,-12)(75,20)
    \path(75,20)(77.8,17.49)(79.02,16.75)(80.06,16.31)(80.98,16.08)(81.83,16)(82.6,16.04)(83.31,16.19)(83.95,16.43)(84.51,16.74)(85.01,17.12)(85.43,17.56)(85.77,18.04)(86.03,18.56)(86.21,19.1)(86.3,19.66)(86.31,20.23)(86.23,20.79)(86.07,21.33)(85.83,21.86)(85.5,22.35)(85.1,22.79)(84.62,23.19)(84.06,23.52)(83.44,23.77)(82.75,23.93)(81.99,24)(81.16,23.95)(80.25,23.75)(79.24,23.36)(78.07,22.69)(75,20)
    \path(75,20)(73,22)
    \put(33,5){\makebox(0,0)[l]{$F_{5}$}}
    \Thicklines
    % (3;1,1)
    \path(5,-30)(5,30)
    \put(6,0){\makebox(0,0)[l]{$-6$}}
    \path(65,-30)(65,30)
    \put(66,0){\makebox(0,0)[l]{$-6$}}
    % (2;1,1)
    \path(15,20)(15,35)
    \put(14.5,31){\makebox(0,0)[r]{$A_{1}$}}
    \path(15,-20)(15,-35)
    \put(14.5,-31){\makebox(0,0)[r]{$A_{2}$}}
    % (3;1,2)
    \path(95,-15)(95,15)
    \path(105,-15)(105,15)
    % [5,-3,5]
    \path(15,-10)(15,10)
    \path(60,-10)(60,10)
    % [13,0,39]
    \path(115,-16)(115,16)
    \path(117,-11)(108,-22)
    \path(110,-17)(110,-32)
    \path(117,11)(108,22)
    \path(110,17)(110,32)
    \path(130,-16)(130,16)
    \path(132,-11)(123,-22)
    \path(125,-17)(125,-32)
    \path(132,11)(123,22)
    \path(125,17)(125,32)
    % (2;1,1)
    \put(65,-45){\arc{15}{3.14}{6.28}}
    \put(65,45){\arc{15}{0}{3.14}}
  \end{picture}
  \caption{$Y_{\I}$}
  \label{fig:Y-I}
\end{figure}

We let $\iota_{\I}$ denote the involution on $Y_{\I}$ induced by the
group~$\Gamma_{\II}$. It is fixed point free. If $N$ is not divisible
by $3$, then $\iota_{\I}$ maps $F_{N}$ to $F_{3N}$ and vice versa. We
denote their common images in $Y_{\II}$ by~$F_{N}$. We have of course
that the self-intersection of $F_{2}$, $F_{13}$ and $F_{5}$
respectively, are the same as for the corresponding curves
on~$Y_{\I}$. The configuration on $Y_{\II}$ is drawn in
figure~\ref{fig:Y-II}.

%%%%%%%%%%%%%%%%%%%%%%%%%%%%%%%%%%%%%%%%%%%%%%%%%%%%%%%%%%%%
% Y_{\II}
\begin{figure}[htbp]
  \centering
  \setlength{\unitlength}{0.9mm}
  \begin{picture}(115,55)(10,-20)
    % F13
    \path(5,25)(130,25)
    \put(5,28){\makebox(0,0)[l]{$F_{13}$}}
    % F2
    \path(15,35)(15,-20)
    \put(16,-20){\makebox(0,0)[l]{$F_{2}$}}
    % F5
    \path(42.38,-20.00)(42.38,12.82)(42.76,21.65)(43.15,25.33)(43.53,27.08)(43.91,28.12)(44.29,28.87)(44.67,29.48)(45.06,30.01)(45.44,30.48)(45.82,30.89)(46.20,31.24)(46.58,31.52)(46.97,31.74)(47.35,31.90)(47.73,31.98)(48.11,32.00)(48.49,31.94)(48.88,31.82)(49.26,31.62)(49.64,31.37)(50.02,31.04)(50.40,30.66)(50.79,30.22)(51.17,29.73)(51.55,29.19)(51.93,28.61)(52.31,27.99)(52.70,27.34)(53.08,26.67)(53.46,25.99)(53.84,25.29)(54.22,24.59)(54.61,23.90)(54.99,23.21)(55.37,22.54)(55.75,21.90)(56.13,21.29)(56.51,20.72)(56.90,20.19)(57.28,19.70)(57.66,19.27)(58.04,18.90)(58.42,18.59)(58.81,18.34)(59.19,18.16)(59.57,18.04)(59.95,18.00)(60.33,18.03)(60.72,18.12)(61.10,18.29)(61.48,18.52)(61.86,18.82)(62.24,19.17)(62.63,19.59)(63.01,20.06)(63.39,20.58)(63.77,21.15)(64.15,21.75)(64.54,22.38)(64.92,23.04)(65.30,23.72)(65.68,24.42)(66.00,25.00)(65.00,26.96)(64.64,27.85)(64.39,28.59)(64.23,29.24)(64.11,29.84)(64.03,30.40)(63.99,30.93)(63.97,31.43)(63.99,31.90)(64.02,32.35)(64.08,32.77)(64.16,33.16)(64.25,33.52)(64.36,33.86)(64.48,34.16)(64.62,34.44)(64.77,34.68)(64.93,34.90)(65.09,35.08)(65.26,35.23)(65.44,35.35)(65.63,35.43)(65.81,35.48)(66.00,35.50)(66.19,35.48)(66.37,35.43)(66.56,35.35)(66.74,35.23)(66.91,35.08)(67.07,34.90)(67.23,34.68)(67.38,34.44)(67.52,34.16)(67.64,33.86)(67.75,33.52)(67.84,33.16)(67.92,32.77)(67.98,32.35)(68.01,31.90)(68.03,31.43)(68.01,30.93)(67.97,30.40)(67.89,29.84)(67.77,29.24)(67.61,28.59)(67.36,27.85)(67.00,26.96)(66.00,25.00)(66.00,25.00)(66.38,24.30)(66.76,23.61)(67.15,22.93)(67.53,22.27)(67.91,21.64)(68.29,21.05)(68.67,20.49)(69.06,19.98)(69.44,19.52)(69.82,19.11)
    \put(41,-20){\makebox(0,0)[r]{$F_{5}$}}
    \Thicklines
    \put(28,29){\arc{15}{0}{3.14}}
    \path(10,0)(30,0)
    \put(20,-3){\makebox(0,0)[l]{$A$}}
    \put(44,0){\arc{15}{1.57}{4.71}}
    \path(44,7.5)(49,7.5)
    \path(44,-7.5)(49,-7.5)
    \put(38,5){\makebox(0,0)[r]{$D$}}

    \put(82,29){\arc{15}{0}{3.14}}
    % (3;1,1)
    \path(55,12)(15,12)
    \put(15,0){\arc{24}{1.57}{4.71}}
    \path(15,-12)(55,-12)
    \put(7,12){\makebox(0,0)[r]{$-6$}}
    \put(110.5,20){\arc{20}{3.14}{6.28}}
    \path(100.5,20)(100.5,15)
    \path(120.5,20)(120.5,15)
    \path(116,20)(130,15)
    \path(127,19)(127,2)
    \path(105,20)(91,15)
    \path(94,19)(94,2)
    % Fixpunkter
    \put(16,2){\makebox(0,0)[l]{$y_{1}$}}
    \put(15,0){\circle*{.8}}
    \put(26,2){\makebox(0,0)[l]{$y_{2}$}}
    \put(25,0){\circle*{.8}}

    \put(47,5){\makebox(0,0)[l]{$y_{3}$}}
    \put(46,7.5){\circle*{.8}}
    \put(47,-5){\makebox(0,0)[l]{$y_{4}$}}
    \put(46,-7.5){\circle*{.8}}

    \put(46,14){\makebox(0,0)[l]{$y_{5}$}}
    \put(45,12){\circle*{.8}}
    \put(46,-14.5){\makebox(0,0)[l]{$y_{6}$}}
    \put(45,-12){\circle*{.8}}
  \end{picture}
  \caption{$Y_{\II}$}
  \label{fig:Y-II}
\end{figure}

The involution on $Y_{\II}$ induced by $\Gamma_{\III}$ has the curve
$F_{13}$ and the $6$ isolated points $y_{1},\dots,y_{6}$ as fixed
point locus. We blow up these points and then we take the quotient.
The surface we get is $Y_{\III}$ blown up 5 times, see
figure~\ref{fig:Y-tilde-III}. On this surface, we have that $F_{2}$,
$F_{13}$ and $F_{5}$ are non-singular rational curves and we have
$F_{2}^{2}=-1$, $F_{13}^{2}=-8$ and $F_{5}^{2}=0$.

%%%%%%%%%%%%%%%%%%%%%%%%%%%%%%%%%%%%%%%%%%%%%%%%%%%%%%%%%%%%
% Y_{\III} uppblåst
\begin{figure}[htbp]
  \centering
  \setlength{\unitlength}{0.9mm}
  \begin{picture}(115,60)(10,-25)
    % F13
    \path(5,25)(130,25)
    \put(5,28){\makebox(0,0)[l]{$F_{13}$}}
    % F2
    \path(15,35)(15,-25)
    \put(16,-25){\makebox(0,0)[l]{$F_{2}$}}
    % F5
    \path(60.0,-25.00)(60.2,-1.79)(60.3,9.44)(60.4,16.37)(60.5,20.70)(60.6,23.43)(60.7,25.20)(60.8,26.37)(60.9,27.17)(61.0,27.74)(61.1,28.17)(61.2,28.50)(61.3,28.77)(61.4,29.00)(61.5,29.19)(61.6,29.36)(61.7,29.51)(61.8,29.64)(61.9,29.76)(62.0,29.86)(62.1,29.94)(62.2,30.02)(62.3,30.08)(62.4,30.13)(62.5,30.17)(62.6,30.20)(62.7,30.21)(62.8,30.22)(62.9,30.22)(63.0,30.20)(63.1,30.18)(63.2,30.15)(63.3,30.11)(63.4,30.06)(63.5,30.00)(63.6,29.94)(63.7,29.87)(63.8,29.79)(63.9,29.70)(64.0,29.61)(64.1,29.51)(64.2,29.40)(64.3,29.29)(64.4,29.18)(64.5,29.05)(64.6,28.93)(64.7,28.80)(64.8,28.66)(64.9,28.52)(65.0,28.37)(65.1,28.23)(65.2,28.08)(65.3,27.92)(65.4,27.76)(65.5,27.60)(65.6,27.44)(65.7,27.27)(65.8,27.11)(65.9,26.94)(66.0,26.76)(66.1,26.59)(66.2,26.42)(66.3,26.24)(66.4,26.07)(66.5,25.89)(66.6,25.71)(66.7,25.53)(66.8,25.36)(66.9,25.18)(67.0,25.00)(67.1,24.82)(67.2,24.65)(67.3,24.47)(67.4,24.30)(67.5,24.12)(67.6,23.95)(67.7,23.78)(67.8,23.61)(67.9,23.44)(68.0,23.27)(68.1,23.11)(68.2,22.94)(68.3,22.78)(68.4,22.63)(68.5,22.47)(68.6,22.32)(68.7,22.17)(68.8,22.02)(68.9,21.87)(69.0,21.73)(69.1,21.59)(69.2,21.46)(69.3,21.33)(69.4,21.20)(69.5,21.07)(69.6,20.95)(69.7,20.83)(69.8,20.72)(69.9,20.61)(70.0,20.50)(70.1,20.40)(70.2,20.30)(70.3,20.20)(70.4,20.11)(70.5,20.02)(70.6,19.94)(70.7,19.86)(70.8,19.79)(70.9,19.72)(71.0,19.65)(71.1,19.59)(71.2,19.54)(71.3,19.48)(71.4,19.44)(71.5,19.39)(71.6,19.35)(71.7,19.32)(71.8,19.29)(71.9,19.26)(72.0,19.24)(72.1,19.22)(72.2,19.21)(72.3,19.20)(72.4,19.20)(72.5,19.20)(72.6,19.20)(72.7,19.21)(72.8,19.23)(72.9,19.24)(73.0,19.27)(73.1,19.29)(73.2,19.32)(73.3,19.36)(73.4,19.40)(73.5,19.44)(73.6,19.49)(73.7,19.54)(73.8,19.59)(73.9,19.65)(74.0,19.71)(74.1,19.77)(74.2,19.84)(74.3,19.91)(74.4,19.99)(74.5,20.07)(74.6,20.15)(74.7,20.23)(74.8,20.32)(74.9,20.41)(75.0,20.50)(75.1,20.60)(75.2,20.69)(75.3,20.79)(75.4,20.89)(75.5,21.00)(75.6,21.10)(75.7,21.21)(75.8,21.32)(75.9,21.43)(76.0,21.54)(76.1,21.66)(76.2,21.77)(76.3,21.89)(76.4,22.00)(76.5,22.12)(76.6,22.24)(76.7,22.35)(76.8,22.47)(76.9,22.59)(77.0,22.71)(77.1,22.82)(77.2,22.94)(77.3,23.05)(77.4,23.17)(77.5,23.28)(77.6,23.39)(77.7,23.50)(77.8,23.60)(77.9,23.71)(78.0,23.81)(78.1,23.91)(78.2,24.01)(78.3,24.10)(78.4,24.19)(78.5,24.28)(78.6,24.37)(78.7,24.45)(78.8,24.52)(78.9,24.59)(79.0,24.66)(79.1,24.72)(79.2,24.78)(79.3,24.83)(79.4,24.87)(79.5,24.91)(79.6,24.94)(79.7,24.97)(79.8,24.98)(79.9,25.00)(80.0,25.00)(80.1,25.00)(80.2,24.98)(80.3,24.96)(80.4,24.93)(80.5,24.90)(80.6,24.85)(80.7,24.79)(80.8,24.72)(80.9,24.65)(81.0,24.56)(81.1,24.46)(81.2,24.35)(81.3,24.23)(81.4,24.09)(81.5,23.95)(81.6,23.79)(81.7,23.62)(81.8,23.43)(81.9,23.23)(82.0,23.02)(82.1,22.79)(82.2,22.55)(82.3,22.29)(82.4,22.02)(82.5,21.73)(82.6,21.43)(82.7,21.10)(82.8,20.76)(82.9,20.41)(83.0,20.03)
    \put(61,-25){\makebox(0,0)[l]{$F_{5}$}}
    \Thicklines
    \put(36,29){\arc{15}{0}{3.14}}
    \put(43,32){\makebox(0,0){$-1$}}
    \path(10,10)(23,18)
    \put(9,10){\makebox(0,0)[r]{$-2$}}
    \path(18,18)(31,10)
    \put(25,15){\makebox(0,0)[l]{$-2$}}
    \path(26,10)(39,18)
    \put(32,13){\makebox(0,0)[l]{$-2$}}
    \path(45,15)(75,15)
    \put(54,17){\makebox(0,0){$-2$}}
    \path(44,5)(49,18)
    \put(47,10){\makebox(0,0)[l]{$-2$}}
    \path(76,5)(71,18)
    \put(73,10){\makebox(0,0)[r]{$-2$}}
    \put(93,29){\arc{15}{0}{3.14}}
    \put(86,32){\makebox(0,0){$-1$}}
    % (3;1,1)
    \path(5,-15)(80,-15)
    \put(28,-18){\makebox(0,0){$-4$}}
    \path(49,-18)(44,-5)
    \put(47,-10){\makebox(0,0)[l]{$-2$}}
    \path(71,-18)(76,-5)
    \put(72,-10){\makebox(0,0)[r]{$-2$}}
    \put(110.5,20){\arc{20}{3.14}{6.28}}
    \path(100.5,20)(100.5,15)
    \path(120.5,20)(120.5,15)
    \put(110,32){\makebox(0,0){$-1$}}
    \path(116,20)(130,15)
    \put(125,14){\makebox(0,0)[r]{$-2$}}
    \path(127,19)(127,2)
    \put(128,10){\makebox(0,0)[l]{$-2$}}
  \end{picture}
  \caption{$\widetilde{Y}_{\III}$}
  \label{fig:Y-tilde-III}
\end{figure}
We blow down the exceptional curve configurations which meet $F_{13}$
and we get the surface~$Y_{\III}$. See figure~\ref{fig:Y-III}.
%%%%%%%%%%%%%%%%%%%%%%%%%%%%%%%%%%%%%%%%%%%%%%%%%%%%%%%%%%%%
% Y_{\III}
\begin{figure}[htbp]
  \centering
  \setlength{\unitlength}{0.9mm}
  \begin{picture}(115,60)(10,-25)
    % F13
    \path(5,25)(26,25)
    \path(26.00,25.00)(28.13,25.04)(30.05,25.16)(31.77,25.35)(33.29,25.62)(34.62,25.95)(35.78,26.36)(36.78,26.81)(37.61,27.32)(38.29,27.87)(38.83,28.45)(39.24,29.06)(39.53,29.69)(39.70,30.31)(39.76,30.94)(39.73,31.55)(39.61,32.13)(39.41,32.68)(39.14,33.19)(38.80,33.64)(38.42,34.05)(37.99,34.38)(37.52,34.65)(37.03,34.84)(36.52,34.96)(36.00,35.00)(35.48,34.96)(34.97,34.84)(34.48,34.65)(34.01,34.38)(33.58,34.05)(33.20,33.64)(32.86,33.19)(32.59,32.68)(32.39,32.13)(32.27,31.55)(32.24,30.94)(32.30,30.31)(32.47,29.69)(32.76,29.06)(33.17,28.45)(33.71,27.87)(34.39,27.32)(35.22,26.81)(36.22,26.36)(37.38,25.95)(38.71,25.62)(40.23,25.35)(41.95,25.16)(43.87,25.04)(46.00,25.00)
    \path(46,25)(86,25)
    \path(86.00,25.00)(88.13,25.04)(90.05,25.16)(91.77,25.35)(93.29,25.62)(94.62,25.95)(95.78,26.36)(96.78,26.81)(97.61,27.32)(98.29,27.87)(98.83,28.45)(99.24,29.06)(99.53,29.69)(99.70,30.31)(99.76,30.94)(99.73,31.55)(99.61,32.13)(99.41,32.68)(99.14,33.19)(98.80,33.64)(98.42,34.05)(97.99,34.38)(97.52,34.65)(97.03,34.84)(96.52,34.96)(96.00,35.00)(95.48,34.96)(94.97,34.84)(94.48,34.65)(94.01,34.38)(93.58,34.05)(93.20,33.64)(92.86,33.19)(92.59,32.68)(92.39,32.13)(92.27,31.55)(92.24,30.94)(92.30,30.31)(92.47,29.69)(92.76,29.06)(93.17,28.45)(93.71,27.87)(94.39,27.32)(95.22,26.81)(96.22,26.36)(97.38,25.95)(98.71,25.62)(100.23,25.35)(101.95,25.16)(103.87,25.04)(106.00,25.00)
    \path(106,25)(123,25)
    \path(122.730,25.000)(123.080,25.000)(123.430,24.990)(123.780,24.980)(124.130,24.950)(124.470,24.910)(124.800,24.850)(125.130,24.760)(125.450,24.650)(125.760,24.500)(126.070,24.320)(126.370,24.110)(126.660,23.870)(126.950,23.580)(127.220,23.260)(127.480,22.900)(127.730,22.510)(127.970,22.080)(128.200,21.610)(128.420,21.120)(128.630,20.590)(128.820,20.040)(129.000,19.470)(129.160,18.880)(129.320,18.270)(129.450,17.660)(129.570,17.030)(129.680,16.400)(129.770,15.770)(129.850,15.150)(129.910,14.540)(129.960,13.930)(129.990,13.350)(130.000,12.780)(130.000,12.240)(129.980,11.720)(129.940,11.230)(129.890,10.780)(129.830,10.350)(129.740, 9.970)(129.650, 9.620)(129.540, 9.310)(129.410, 9.040)(129.270, 8.810)(129.110, 8.630)(128.940, 8.490)(128.760, 8.390)(128.560, 8.340)(128.350, 8.340)(128.130, 8.380)(127.890, 8.460)(127.650, 8.590)(127.390, 8.770)(127.130, 8.990)(126.850, 9.250)(126.570, 9.550)(126.270, 9.890)(125.970,10.270)(125.660,10.690)(125.340,11.140)(125.020,11.620)(124.690,12.130)(124.350,12.670)(124.010,13.230)(123.670,13.820)(123.320,14.410)(122.960,15.030)(122.610,15.650)(122.250,16.280)(121.890,16.900)(121.520,17.530)(121.160,18.150)(120.790,18.760)(120.430,19.350)(120.060,19.930)(119.690,20.490)(119.320,21.010)(118.960,21.520)(118.590,21.990)(118.230,22.420)(117.860,22.830)(117.500,23.190)(117.140,23.520)(116.790,23.810)(116.430,24.070)(116.080,24.280)(115.730,24.470)(115.380,24.620)(115.040,24.740)(114.700,24.830)(114.360,24.900)(114.020,24.950)(113.690,24.980)(113.360,24.990)(113.040,25.000)(112.720,25.000)(112.400,25.000)(112.090,25.000)(111.780,24.980)(111.480,24.960)(111.180,24.920)(110.880,24.860)(110.590,24.780)(110.300,24.670)(110.010,24.530)(109.730,24.360)(109.450,24.160)(109.180,23.920)(108.910,23.640)(108.650,23.330)(108.390,22.980)(108.130,22.590)(107.880,22.160)(107.630,21.710)(107.380,21.220)(107.140,20.700)(106.900,20.160)(106.670,19.590)(106.440,19.000)(106.210,18.400)(105.990,17.780)
    \put(5,28){\makebox(0,0)[l]{$F_{13}$}}
    % F2
    \path(15,35)(15,-25)
    \put(16,-25){\makebox(0,0)[l]{$F_{2}$}}
    % F5
    \path(60.0,-25.00)(60.2,-1.79)(60.3,9.44)(60.4,16.37)(60.5,20.70)(60.6,23.43)(60.7,25.20)(60.8,26.37)(60.9,27.17)(61.0,27.74)(61.1,28.17)(61.2,28.50)(61.3,28.77)(61.4,29.00)(61.5,29.19)(61.6,29.36)(61.7,29.51)(61.8,29.64)(61.9,29.76)(62.0,29.86)(62.1,29.94)(62.2,30.02)(62.3,30.08)(62.4,30.13)(62.5,30.17)(62.6,30.20)(62.7,30.21)(62.8,30.22)(62.9,30.22)(63.0,30.20)(63.1,30.18)(63.2,30.15)(63.3,30.11)(63.4,30.06)(63.5,30.00)(63.6,29.94)(63.7,29.87)(63.8,29.79)(63.9,29.70)(64.0,29.61)(64.1,29.51)(64.2,29.40)(64.3,29.29)(64.4,29.18)(64.5,29.05)(64.6,28.93)(64.7,28.80)(64.8,28.66)(64.9,28.52)(65.0,28.37)(65.1,28.23)(65.2,28.08)(65.3,27.92)(65.4,27.76)(65.5,27.60)(65.6,27.44)(65.7,27.27)(65.8,27.11)(65.9,26.94)(66.0,26.76)(66.1,26.59)(66.2,26.42)(66.3,26.24)(66.4,26.07)(66.5,25.89)(66.6,25.71)(66.7,25.53)(66.8,25.36)(66.9,25.18)(67.0,25.00)(67.1,24.82)(67.2,24.65)(67.3,24.47)(67.4,24.30)(67.5,24.12)(67.6,23.95)(67.7,23.78)(67.8,23.61)(67.9,23.44)(68.0,23.27)(68.1,23.11)(68.2,22.94)(68.3,22.78)(68.4,22.63)(68.5,22.47)(68.6,22.32)(68.7,22.17)(68.8,22.02)(68.9,21.87)(69.0,21.73)(69.1,21.59)(69.2,21.46)(69.3,21.33)(69.4,21.20)(69.5,21.07)(69.6,20.95)(69.7,20.83)(69.8,20.72)(69.9,20.61)(70.0,20.50)(70.1,20.40)(70.2,20.30)(70.3,20.20)(70.4,20.11)(70.5,20.02)(70.6,19.94)(70.7,19.86)(70.8,19.79)(70.9,19.72)(71.0,19.65)(71.1,19.59)(71.2,19.54)(71.3,19.48)(71.4,19.44)(71.5,19.39)(71.6,19.35)(71.7,19.32)(71.8,19.29)(71.9,19.26)(72.0,19.24)(72.1,19.22)(72.2,19.21)(72.3,19.20)(72.4,19.20)(72.5,19.20)(72.6,19.20)(72.7,19.21)(72.8,19.23)(72.9,19.24)(73.0,19.27)(73.1,19.29)(73.2,19.32)(73.3,19.36)(73.4,19.40)(73.5,19.44)(73.6,19.49)(73.7,19.54)(73.8,19.59)(73.9,19.65)(74.0,19.71)(74.1,19.77)(74.2,19.84)(74.3,19.91)(74.4,19.99)(74.5,20.07)(74.6,20.15)(74.7,20.23)(74.8,20.32)(74.9,20.41)(75.0,20.50)(75.1,20.60)(75.2,20.69)(75.3,20.79)(75.4,20.89)(75.5,21.00)(75.6,21.10)(75.7,21.21)(75.8,21.32)(75.9,21.43)(76.0,21.54)(76.1,21.66)(76.2,21.77)(76.3,21.89)(76.4,22.00)(76.5,22.12)(76.6,22.24)(76.7,22.35)(76.8,22.47)(76.9,22.59)(77.0,22.71)(77.1,22.82)(77.2,22.94)(77.3,23.05)(77.4,23.17)(77.5,23.28)(77.6,23.39)(77.7,23.50)(77.8,23.60)(77.9,23.71)(78.0,23.81)(78.1,23.91)(78.2,24.01)(78.3,24.10)(78.4,24.19)(78.5,24.28)(78.6,24.37)(78.7,24.45)(78.8,24.52)(78.9,24.59)(79.0,24.66)(79.1,24.72)(79.2,24.78)(79.3,24.83)(79.4,24.87)(79.5,24.91)(79.6,24.94)(79.7,24.97)(79.8,24.98)(79.9,25.00)(80.0,25.00)(80.1,25.00)(80.2,24.98)(80.3,24.96)(80.4,24.93)(80.5,24.90)(80.6,24.85)(80.7,24.79)(80.8,24.72)(80.9,24.65)(81.0,24.56)(81.1,24.46)(81.2,24.35)(81.3,24.23)(81.4,24.09)(81.5,23.95)(81.6,23.79)(81.7,23.62)(81.8,23.43)(81.9,23.23)(82.0,23.02)(82.1,22.79)(82.2,22.55)(82.3,22.29)(82.4,22.02)(82.5,21.73)(82.6,21.43)(82.7,21.10)(82.8,20.76)(82.9,20.41)(83.0,20.03)
    \put(61,-25){\makebox(0,0)[l]{$F_{5}$}}
    \Thicklines
    \path(10,10)(23,18)
    \put(9,10){\makebox(0,0)[r]{$-2$}}
    \path(18,18)(31,10)
    \put(25,15){\makebox(0,0)[l]{$-2$}}
    \path(26,10)(39,18)
    \put(32,13){\makebox(0,0)[l]{$-2$}}
    \path(45,15)(75,15)
    \put(54,17){\makebox(0,0){$-2$}}
    \path(44,5)(49,18)
    \put(47,10){\makebox(0,0)[l]{$-2$}}
    \path(76,5)(71,18)
    \put(73,10){\makebox(0,0)[r]{$-2$}}
    % (3;1,1)
    \path(5,-15)(80,-15)
    \put(28,-18){\makebox(0,0){$-4$}}
    \path(49,-18)(44,-5)
    \put(47,-10){\makebox(0,0)[l]{$-2$}}
    \path(71,-18)(76,-5)
    \put(72,-10){\makebox(0,0)[r]{$-2$}}
  \end{picture}
  \caption{$Y_{\III}$}
  \label{fig:Y-III}
\end{figure}

Since $\zeta_{k}(-1)= \frac{1}{6}$, we get that
$\int_{X}\omega=\frac{4}{3}$ by equation~\eqref{eq:invariants-volume}.
Furthermore $e_{2}=8$ and $e_{3}=4$, so
by~\eqref{eq:invariants-Euler}, we get $e(X)=8$. It is now
straightforward to compute the numerical invariants of all surfaces.
We summarise the result in the following table:
\begin{center}
  \renewcommand{\arraystretch}{1.1}
  \begin{tabular}{c|cccc}
    & $Y$ & $Y_{\I}$ & $Y_{\II}$ & $Y_{\III}$ \\
    \hline
    $e$ & $22$ & $28$ & $14$ & $12$ \\
    $\chi$ & $2$ & $2$ & $1$ & $1$ \\
    $K^{2}$ & $2$ & $-4$ & $-2$ & $0$ \\
    $q$ & $0$ & $0$ & $0$ & $0$ \\
    $p_{g}$ & $1$ & $1$ & $0$ & $0$ \\
  \end{tabular}
\end{center}

We have the following corollary of Castelnuovo's criterion for
rationality of a surface (see~\cite{hirzebruch1}, corollary~1,
p.~255):
\begin{proposition}
  \label{prop:castelnouvos-criterion}
  Let $Z$ be a non-singular surface with $q=0$. If $Z$ contains a
  non-singular rational curve $D$ with $D^{2}\geq0$, then $Z$ is
  rational.
\end{proposition}
This directly applies to give that $Y_{\III}$ is a rational surface.
Namely, if we blow down the connected configuration consisting of
$F_{2}$ and three $(-2)$-curves, then the image of the $(-4)$-curve is
a non-singular rational curve with self-intersection~$0$.

Studying the intersection matrix of the divisors generated by the
curves in figure~\ref{fig:Y-tilde-III}, we find that these curves
generate a sublattice of $\Pic(\widetilde{Y}_{\III})$ of index~2.
Using the adjunction formula, one gets the following representation of
the canonical divisor of $\widetilde{Y}_{\III}$:
\begin{equation}
  \label{eq:2K}
  2K=4F_{2}-F_{13}+3C_{14}+2C_{15}+C_{16}.
\end{equation}
Using~\eqref{eq:2K}, we can trace $2K$ backwards and get the
following relation on~$Y_{\II}$:
\begin{equation}
  \label{eq:2K-on-Y-II}
  2K=4F_{2}+2A.
\end{equation}
Let now $Z$ denote the surface we get if we blow down first $F_{2}$
and then $A$ on~$Y_{\II}$. By~\eqref{eq:2K-on-Y-II}, we have that
$2K=0$ on~$Z$. Furthermore, $e(Z)=12$ and $q(Z)=0$. This shows that
$Z$ is an Enriques surface (see~\cite{barth-peters-ven}, chapter~VI).
Since the map $Y_{\I}\rightarrow Y_{\II}$ is unramified, we
immediately get that $Y_{\I}$ is a $K3$-surface blown up $4$ times
with
\begin{equation}
  \label{eq:K-on-Y-I}
  K=2F_{2}+A_{1}+2F_{6}+A_{2}.
\end{equation}
This implies that the canonical divisor on $Y$ is given by
\begin{displaymath}
  K=2F_{2}+2F_{6}+\sum A_{i}+\sum B_{i}.
\end{displaymath}
We conclude that $Y$ is a minimal surface (since if $E$  an exceptional
divisor on $Y$, then $KE=-1$, which is impossible since $K$ is
effective and does not contain any exceptional component). Furthermore,
since $K$ is effective and $K^{2}>0$, we have that $Y$ is of general
type (see~\cite{barth-peters-ven}).

%%%%%%%%%%%%%%%%%%%%%%%%%%%%%%%%%%%%%%%%%%%%%%%%%%%%%%%%%%%%%%%%%%%%%%

\end{document}